\documentclass[11pt]{article}
\usepackage[T1]{fontenc} 
\usepackage[utf8]{inputenc}

\usepackage{lmodern}

\usepackage{secdot}
\usepackage{amssymb}
\usepackage{amsmath}
\usepackage{latexsym}
\usepackage{calrsfs}
\usepackage{graphicx,color,psfrag}
\usepackage{natbib}

\newtheorem{theorem}{Theorem}[section]
\newtheorem{lemma}[theorem]{Lemma}

\newtheorem{proposition}[theorem]{Proposition}

\newenvironment{proof}{\trivlist\item[]\emph{Proof}.}
{\unskip\nobreak\hskip 1em plus 1fil\nobreak$\Box$
\parfillskip=0pt\endtrivlist}

\def\argmax{\mathop{\rm arg\,max}}

\newcommand{\Ex}{\ensuremath{\mathsf{E}}}
\newcommand{\transp}[1]{{#1}^{\scriptscriptstyle{\mathsf{T}}}}

\usepackage[colorlinks=true,urlcolor=blue,citecolor=blue,linkcolor=blue]{hyperref}
\usepackage[noabbrev,nameinlink]{cleveref}     

\def\EMAIL#1{\href{mailto:#1}{#1}}

\begin{document}

\title{A {\Large $(2/3)n^3$} fast-pivoting algorithm for the Gittins index and
  optimal stopping of a Markov chain}


\author{Jos\'e Ni\~no-Mora 
\\ Department of Statistics \\
    Carlos III University of Madrid \\
     28903 Getafe (Madrid), Spain \\  \EMAIL{jose.nino@uc3m.es}, \href{http://alum.mit.edu/www/jnimora}{http://alum.mit.edu/www/jnimora} \\
      ORCID: \href{http://orcid.org/0000-0002-2172-3983}{0000-0002-2172-3983}}
 
\date{Published in \textit{INFORMS Journal on Computing}, vol.\  19, pp. 596--606,  2007 \\ \vspace{.1in}
DOI: \href{https://doi.org/10.1287/ijoc.1060.0206}{10.1287/ijoc.1060.0206}}

\maketitle


\begin{abstract}%
This paper presents a new \emph{fast-pivoting}  algorithm that computes the $n$ Gittins index values of
 an $n$-state bandit --- in the discounted and undiscounted cases ---
by performing
$(2/3) n^3 + O(n^2)$ arithmetic operations, thus attaining better
 complexity than
 previous algorithms and  matching 
 that of  solving a corresponding linear-equation system by Gaussian
elimination.
The algorithm further applies to the problem of optimal stopping of a
 Markov chain, for which a novel Gittins-index solution approach is introduced.
The algorithm draws on Gittins and Jones' (1974) index definition via 
calibration, on Kallenberg's (1986) proposal of using parametric
 linear programming, on Dantzig's simplex method, on Varaiya et al.'s (1985)
 algorithm,  and on the author's earlier work.
The paper elucidates the structure of parametric simplex
 tableaux. Special structure is exploited to reduce
 the computational effort of pivot steps,
decreasing the operation count by a factor of three relative to 
using conventional pivoting, and by a factor of $3/2$ relative to
recent state-elimination algorithms.
A computational study demonstrates 
 significant time savings against alternative algorithms.
\end{abstract}%


\newpage

\section{Introduction}
\label{s:intro}
Consider the following \emph{optimal-stopping problem}. 
A Markov chain with state $X(t)$ evolves through the finite or countable
state space $N$, according to transition
probabilities $p_{ij}$.
At each discrete time period $t \geq 0$, one must decide either to let
the chain continue or to stop it, 
based on information so far.
Continuing the chain when it occupies state $i$ yields an immediate
reward $R_i$, but a charge $\nu$ is incurred. 
If the chain is instead stopped, a terminal reward $Q_i$ is earned.
We assume the $R_i$'s and $Q_i$'s to
be uniformly bounded in the countable-state case.
Rewards earned over time
 are discounted with factor $0 < \beta \leq 1$.
We want to find an optimal
 \emph{stopping rule} that maximizes the expected total discounted
value of rewards earned minus charges incurred. We can formulate
such problem as 
\begin{equation}
\label{eq:prob1}
\max_{\tau \geq 0} \Ex_i^\tau\left[\sum_{t=0}^{\tau-1} (R_{X(t)} - \nu)
  \beta^t + Q_{X(\tau)} \beta^\tau\right],
\end{equation}
where $\Ex_i^\tau[\cdot]$ denotes
expectation starting at $i$ and, abusing notation, $\tau$ denotes both
a stopping rule  and the corresponding 
\emph{stopping time} $0 \leq \tau \leq +\infty$.

An important application area for optimal stopping
 is computational finance, where several fundamental
problems are formulated as special cases of (\ref{eq:prob1}),
 such as
deciding the optimal time to sell a stock, or  to exercise a perpetual
American option.
While prevailing 
models  are continuous, there is emerging interest in 
more numerically tractable Markov-chain
approximations, to which the results in this paper are applicable.
See \citet{duangausim03}.

Consider now the special zero-terminal-rewards case of problem (\ref{eq:prob1}), i.e., $Q_j \equiv 0$:
\begin{equation}
\label{eq:prob2}
\max_{\tau \geq 0} \Ex_i^\tau\left[\sum_{t=0}^{\tau-1} (R_{X(t)} - \nu) \beta^t\right].
\end{equation}
 \citet{gijo74}  showed --- in the reformulation of
problem (\ref{eq:prob2}) where $\nu$ plays the role of a subsidy per
passive period --- that one can attach to every state
$i$ a quantity $\nu_i^*$, now termed the \emph{Gittins index}, 
such  that it is optimal to stop 
at state $i$ if and only if $\nu_i^* \leq \nu$.
Such an index result was first obtained by
\citet{bjk56} in the classical setting of a Bayesian Bernoulli bandit.
\citet{whit80} highlighted and exploited the alternative interpretation of the 
Gittins index as a critical \emph{constant} terminal reward in an optimal-stopping
problem --- think of the reformulation of (\ref{eq:prob2})
where a terminal reward of $\nu/(1-\beta)$ is earned.

\citet{gijo74}  further proved that such an index furnishes an
efficient solution to the \emph{multiarmed-bandit problem}, 
concerning the optimal dynamic effort allocation to a collection of 
bandits, 
one of which must be engaged at a time: one should always engage
a bandit with largest index.
\citet{gi79} showed that
\begin{equation}
\label{eq:girep}
\nu_i^* = \max_{\tau > 0} \frac{\displaystyle \Ex_i^\tau\left[\sum_{t=0}^{\tau-1}
    R_{X(t)} \beta^t\right]}{\displaystyle \Ex_i^\tau\left[\sum_{t=0}^{\tau-1} \beta^t\right]},
\end{equation}
so that $\nu_i^*$ represents
 the  maximum rate of expected discounted reward per unit of
expected discounted time that can be achieved starting at $i$.
The corresponding \emph{undiscounted Gittins index} is obtained by
setting $\beta = 1$ above.

This paper presents two main contributions: (i) we show that 
problem (\ref{eq:prob1}) can be reduced to (\ref{eq:prob2}),
and can thus be solved via the Gittins index; and (ii) 
for a bandit with a finite number $n$ of states, we give an 
efficient algorithm based on the parametric version of Dantzig's 
simplex method
 that computes the Gittins index by
performing 
 $(2/3)n^3 + O(n^2)$ arithmetic operations, thus achieving better
complexity than previous algorithms.
It appears unlikely that such complexity can be improved, as it
matches that of solving a corresponding linear-equation system by 
Gaussian elimination.

We next review previous algorithmic approaches to 
 (\ref{eq:prob1}) and (\ref{eq:prob2}). 
The optimal stopping of a Markov chain 
is a fundamental problem that has
been the subject of extensive research attention. See, e.g., 
\citet{sonin99} and the references therein.
Conventional  approaches are based on formulating the
Bellman equations, which are then solved by the
value-iteration, policy-iteration or linear-programming (LP) 
methods, where the required number of
iterations is not clear a priori.
In contrast,  \citet{sonin99}
introduced a one-pass algorithm based on state-elimination ideas,
which recursively solves an $n$-state problem in $n$ iterations.

Regarding the index,  \citet{vawabu} introduced an $n$-step
algorithm that performs $(1/3) n^4 + O(n^3)$ arithmetic operations.
Both this and the  \citet{kl} index algorithm were
elucidated in \citet{beni} as special cases of
an \emph{adaptive-greedy algorithm} that solves the underlying polyhedral 
LP formulation.

 \citet{chenka86} showed that the 
index of a \emph{fixed state} can be computed by 
solving an LP problem with $n+1$ variables and 
$n$ constraints. 
 \citet{kallenb86} proposed computing the index by
solving a parametric LP problem as in \citet{saatgass54}. 
Such an approach involves $n$ pivot steps of the parametric version of
Dantzig's
simplex method, which requires
$2 n^3 + O(n^2)$ arithmetic operations.

More recently, \citet{kattaset04} and \citet{sonin05}
introduced index algorithms based on 
state elimination, which perform
$n^3 + O(n^2)$ operations.

The new algorithm we develop in this paper 
draws on  \citet{gi79},  \citet{kallenb86}, Dantzig's simplex method, 
\citet{vawabu},  and 
  \citet{nmaap01,nmmp02,nmmor06}.
We elucidate the structure of the parametric simplex
 tableau. By exploiting special structure,
 the computational effort of pivot steps is reduced,
decreasing the operation count by a factor of three relative to 
using conventional pivoting, and by a factor of $3/2$ relative to
state-elimination algorithms.
Such a \emph{fast-pivoting} algorithm applies both to the discounted and the undiscounted
Gittins index.
A computational study demonstrates that  the algorithm achieves in practice
 significant time savings against conventional-pivoting and
 state-elimination algorithms.

We further give a variant of the algorithm that efficiently
computes additional quantities that are useful for other purposes. 
We have used such a variant in
\citet{nmjoc06} as the first stage of a two-stage method to 
compute the 
\emph{switching index} of
 \citet{asatene96} efficiently for bandits with switching costs.

Section~\ref{s:sospgi} shows that optimal stopping problem
(\ref{eq:prob1}) can be solved via the Gittins index.
Section~\ref{s:mgika} describes the bandit model of concern, 
reviews a convenient definition of the Gittins index, and discusses
the algorithm of \citet{vawabu}.
Section~\ref{s:ecge} elucidates the structure of the LP parametric 
simplex tableau for computing the index and formulates the 
conventional-pivoting algorithm. 
Section~\ref{s:ess} exploits special structure to
reduce the computational effort of pivot steps.
Section~\ref{s:piat} describes the new fast-pivoting algorithm.
Section~\ref{s:rsea} discusses relations with the 
state-elimination algorithm. 
Section~\ref{s:cargi} reviews the concept of undiscounted Gittins
index, and shows how to compute it with the fast-pivoting algorithm.
Section~\ref{s:ce} reports the results of a computational
study.
Section~\ref{sec-Conclusions} concludes.

\section{Gittins-Index Solution to Optimal Stopping}
\label{s:sospgi}
In this section we assume that the chain's state space $N$ is either
   finite or countably infinite and, in the latter case,
 that rewards $R_i$ and $Q_i$ are bounded.
Consider optimal-stopping problem (\ref{eq:prob1}),
in the discounted case $0 < \beta < 1$.
The discounted value of rewards earned under stopping rule
 $\tau$ starting at $i$ is
given by the \emph{reward measure}
\begin{equation}
\label{eq:fipi}
f_i^\tau(\mathbf{R}, \mathbf{Q})
\triangleq \Ex_i^\tau\left[\sum_{t=0}^{\tau-1} R_{X(t)} \beta^t +
  Q_{X(\tau)} \beta^\tau\right],
\end{equation}
where we have made explicit the dependence on continuation rewards $\mathbf{R} = (R_j)_{j
  \in N}$ and terminal rewards $\mathbf{Q} = (Q_j)_{j
  \in N}$. 
The corresponding discounted amount of effort expended is given by the
\emph{work measure}
\begin{equation}
\label{eq:gipi}
g_i^\tau
\triangleq \Ex_i^\tau\left[\sum_{t=0}^{\tau-1} \beta^t\right].
\end{equation}

We can thus restate optimal-stopping problem (\ref{eq:prob1}) as 
\begin{equation}
\label{eq:prob1b}
\max_{\tau \geq 0} f_i^\tau(\mathbf{R}, \mathbf{Q}) - \nu g_i^\tau.
\end{equation}

We now make the key observation that the above discounted 
optimal-stopping problem is readily reformulated as an
  \emph{infinite-horizon Markov
  decision process} (MDP) with two actions available at each state and
  time period: $a =
  1$ (active) and $a = 0$ (passive). In the latter,
if the active action is taken in state $i$, reward $R_i$ accrues and the state moves
  according to transition probabilities $p_{ij}$;
if the passive action is taken instead, reward $(1-\beta) Q_i$
  accrues and the state does not change. 
Notice that the scale factor $1-\beta$ incorporated into passive rewards
  accounts for
  the fact that receiving a lump terminal reward of $Q_i$ is
  equivalent to receiving a discounted pension of $(1-\beta) Q_i$ at every period
over an infinite horizon.

The reader might be concerned at this point that such an infinite-horizon
MDP formulation --- which differs from the conventional one based on 
introducing a terminal state ---
allows policies that take the active action after the passive one has
been used, which is not allowed in 
 (\ref{eq:prob1}). 
Such apparent discrepancy is, however, resolved through MDP theory,
 which, since
we assumed bounded rewards,
ensures existence of an \emph{optimal} policy for such an
infinite-horizon MDP that is Markov stationary, and hence consistent with
a stopping rule (including the possibility of 
never stopping).

We next further
 deploy MDP theory in the reformulated infinite-horizon problem.  
Consider the
 \emph{discounted state-action occupancy measures}:
for any stopping rule $\tau$ (now viewed as an MDP
policy in the reformulated model), initial
state $i$, 
action $a$, and state $j$, let
\[
x_{ij}^{a, \tau} \triangleq \Ex_i^\tau\left[\sum_{t=0}^\infty 
 1_{\{a(t) = a, X(t) = j\}} \beta^t\right]
\]
be the expected total discounted 
time expended taking action $a$ in state $j$ under policy $\tau$ starting
at $i$.
It is well known that, for fixed $i$,  the 
$2 n$ measures $x_{ij}^{a, \tau}$, for $(j, a) \in N \times \{0, 1\}$,
satisfy a system of $n$ linear equations.
Writing $\mathbf{x}_i^{a, \tau} = (x_{ij}^{a, \tau})_{j \in N}$ as 
a \emph{row vector}, denoting by $\mathbf{I}$ the identity matrix,
and by $\mathbf{e}_i \in \mathbb{R}^n$ the unit coordinate
row vector having the one in the position of state $i$, such an equation
system
is
\begin{equation}
\label{eq:les}
(1-\beta) \mathbf{x}_i^{0, \tau} + \mathbf{x}_i^{1, \tau} (\mathbf{I} - \beta \mathbf{P}) = \mathbf{e}_i.
\end{equation}

We can formulate reward and work measures as linear functions of occupancies:
\begin{equation}
\label{eq:figilf}
\begin{split}
f_i^\tau(\mathbf{R}, \mathbf{Q}) & = \sum_{j \in N} R_j x_{ij}^{1, \tau} + 
\sum_{j \in N} (1-\beta) Q_j x_{ij}^{0, \tau} = 
\mathbf{x}_i^{1, \tau} \mathbf{R} + 
 (1-\beta) \mathbf{x}_i^{0, \tau} \mathbf{Q}, \\
g_i^\tau & = \sum_{j \in N} x_{ij}^{1, \tau} = \mathbf{x}_i^{1, \tau}
  \mathbf{1}.
\end{split}
\end{equation}

We next draw on the above to show that terminal rewards can be
eliminated from  (\ref{eq:prob1}), both in
the discounted and the undiscounted cases. 

\begin{lemma}
\label{lma:osgi} Under any stopping rule $\tau$ it holds that, for $0
< \beta \leq 1$\textup{:}
\begin{equation}
\label{eq:fiosgi}
f_i^\tau(\mathbf{R}, \mathbf{Q}) = Q_i + f_i^\tau(\mathbf{R} - (\mathbf{I} - \beta \mathbf{P}) \mathbf{Q}, \mathbf{0}).
\end{equation}
\end{lemma}
\begin{proof}
In the discounted case $\beta < 1$, we use the above to write
\begin{align*}
f_i^\tau(\mathbf{R}, \mathbf{Q}) & = \mathbf{x}_i^{1, \tau} \mathbf{R}
+ (1-\beta) \mathbf{x}_i^{0, \tau} \mathbf{Q} \\
& = \mathbf{x}_i^{1, \tau} \mathbf{R} + \{\mathbf{e}_i - 
 \mathbf{x}_i^{1, \tau} (\mathbf{I} - \beta \mathbf{P})\} \mathbf{Q}
 =  Q_i + \mathbf{x}_i^{1, \tau}  \{\mathbf{R} - (\mathbf{I} - \beta
\mathbf{P})\} \mathbf{Q}\} \\
& = Q_i + f_i^\tau(\mathbf{R} -
(\mathbf{I} - \beta \mathbf{P}) \mathbf{Q}, \mathbf{0}).
\end{align*}

Taking the limit as $\beta \nearrow 1$ in the latter identity yields
the result for $\beta = 1$.
\end{proof}

Consider now the bandit with zero passive rewards and
modified active-reward vector
\begin{equation}
\label{eq:rhat}
\hat{\mathbf{R}} \triangleq \mathbf{R} - (\mathbf{I} - \beta \mathbf{P})
\mathbf{Q},
\end{equation}
and let $\hat{\nu}_i^*$ be its Gittins index. 

We next give the main result of this section.
\begin{theorem}
\label{the:osgi}
It is optimal to stop 
at state $i$ in problem \textup{(\ref{eq:prob1})} iff
$\hat{\nu}_i^* \leq \nu$. The minimum optimal-stopping time is thus
$
\tau^* = \min \left\{t \geq 0\colon \hat{\nu}_{X(t)}^* \leq \nu\right\}.
$
\end{theorem}
\begin{proof}
The above discussion shows that  (\ref{eq:prob1}) is
equivalent to a corresponding problem without terminal rewards of the
form (\ref{eq:prob2}), with rewards $\hat{R}_i$ given by
(\ref{eq:rhat}).
Since the latter is solved by the stated Gittins-index 
policy, the result follows.
\end{proof}

\section{Gittins Index and VWB  Algorithm}
\label{s:mgika}
In the remainder of the paper, we focus on a bandit having a finite
number $n$ of states.
This section reviews a convenient definition of the Gittins index, the
concepts of marginal work, reward, and 
productivity measures introduced in \citet{nmaap01,nmmp02,nmmor06},
and the index algorithm of
\citet{vawabu}, to which we will refer as VWB.

In light of Section \ref{s:sospgi}, we need only focus 
 on  (\ref{eq:prob2}), which we formulate as
\begin{equation}
\label{eq:nuwp}
\max_{\tau \geq 0} f_i^\tau - \nu g_i^\tau.
\end{equation}

By MDP theory we know that, for 
any continuation charge $\nu$, 
there exists an optimal stopping rule for
 (\ref{eq:nuwp}) that is Markov-stationary and  independent
of the initial state.
We represent such rules by their \emph{continuation sets}. 
Hence, to any charge $\nu$
 there corresponds a \emph{minimal optimal continuation set}
$S^*(\nu) \subseteq N$. The latter
 are characterized by the \emph{Gittins index}
$\nu_j^*$ that is attached to every state $j$, so that
\[
S^*(\nu) = \left\{j \in N\colon \nu_j^* > \nu\right\}, \quad \nu
\in \mathbb{R}.
\]

For an action $a$ and continuation set $S$, 
let $\langle a, S\rangle$ be the policy that
takes action $a$ in the first period and adopts the 
\emph{$S$-active policy} (which continues over $S$ and 
stops over $S^c \triangleq N \setminus S$) thereafter. 
Define, for a state $i$, 
 the  \emph{$(i, S)$-marginal work measure} 
\begin{equation}
\label{eq:mwm}
w_i^S \triangleq g_i^{\langle 1, S\rangle} - g_i^{\langle 0, S\rangle},
\end{equation}
as the marginal increase in work expended that results from
continuing, instead of stopping in the first period,
starting at $i$, provided that the $S$-active policy is adopted thereafter.
Define further the \emph{$(i, S)$-marginal reward measure} 
\begin{equation}
\label{eq:mrm}
r_i^S \triangleq f_i^{\langle 1, S\rangle} - f_i^{\langle 0, S\rangle}
\end{equation}
as the corresponding marginal increase in rewards earned.
We will see that marginal workloads are positive, which
allows us to 
define the \emph{$(i, S)$-marginal productivity rate} 
\begin{equation}
\label{eq:mpm}
\nu_i^S \triangleq \frac{r_i^S}{w_i^S}.
\end{equation}

We next use the above quantities --- which are readily computed by
solving
linearequation systems, as will be shown in the next section ---
to formulate the version of the VWB algorithm 
in Table \ref{tab:kvwb}.
The algorithm
proceeds in $n$ steps by generating a sequence of states $i_k$ with 
nonincreasing index values and
a corresponding nested sequence of 
continuation sets $S_0 = \emptyset, S_k = \{i_1, \ldots, i_k\}$, for 
$k = 1, \ldots, n$.

\begin{table}[tb]
\caption{Formulation of the VWB Gittins-Index Algorithm via the 
 $w_i^S$'s and $r_i^S$'s.}
\begin{center}
\fbox{%
\begin{minipage}{\textwidth}
{\bf ALGORITHM VWB:} 
\begin{tabbing}
{\bf set }
$S_0 := \emptyset$  \\
{\bf for} \= $k := 1$ {\bf to} $n$ {\bf do} \\
 \> {\bf compute} $w_i^{S_{k-1}}, r_i^{S_{k-1}}, \nu_i^{S_{k-1}} = r_i^{S_{k-1}}/w_i^{S_{k-1}}, \, i
 \in S_{k-1}^c$ \\
 \> {\bf pick} 
 $i_{k} \in \argmax
      \left\{\nu^{S_{k-1}}_i\colon
                i \in S_{k-1}^c\right\}$
  \\
 \> $\nu_{i_k}^* := 
 \nu^{S_{k-1}}_{i_k}, \, \, S_{k} := S_{k-1} \cup \{i_{k}\}$  \\
{\bf end} 
\end{tabbing}
\end{minipage}}
\end{center}
\label{tab:kvwb}
\end{table}

In \citet{vawabu}, the ratios $\nu_i^{S_{k-1}}$ in Table 
\ref{tab:kvwb}
are calculated as $\nu_i^{S_{k-1}} = a^{(k)}_i/b^{(k)}_i$, where
 $a^{(k)}_i$ and $b^{(k)}_i$ are obtained 
by solving the following linearequation systems: for $i \in N$,
\begin{align*}
a_i^{(k)} & = \beta R_i + \beta \sum_{j \in S_{k-1}} p_{ij} a_j^{(k)} \\
b_i^{(k)} & = \beta + \beta \sum_{j \in S_{k-1}} p_{ij} b_j^{(k)}.
\end{align*}

Using the results of the next section and a bit of algebra yields that such quantities are related to 
our $w_i^{S_{k-1}}$ and $r_i^{S_{k-1}}$ above by 
\[
a_i^{(k)} = 
\begin{cases}
\beta r_i^{S_{k-1}} / (1-\beta) & \text{if } i \in S_{k-1} \\
\beta r_i^{S_{k-1}} & \text{if } i \in S_{k-1}^c
\end{cases}
\quad \text{ and } \quad
b_i^{(k)} = 
\begin{cases}
\beta w_i^{S_{k-1}} / (1-\beta) & \text{if } i \in S_{k-1} \\
\beta w_i^{S_{k-1}} & \text{if } i \in S_{k-1}^c.
\end{cases}
\]

We next assess the computational complexity of algorithm VWB.
\begin{proposition}
\label{pro:ccvwb}
Algorithm VWB performs $(1/3) n^4 + O(n^3)$ arithmetic operations.
\end{proposition}
\begin{proof}
The count is dominated by the solution of two linear
equation systems of size $k$ at step $k$, each taking
$(2/3) k^3 + O(k^2)$ arithmetic operations, yielding a total of
\[
2 \sum_{k=2}^{n} \left\{(2/3) k^3 + O(k^2)\right\} = (1/3) n^4 + O(n^3).
\]
\end{proof}

\section{LP Formulation and Parametric Simplex Tableau}
\label{s:ecge}
We set out in this section to formulate 
(\ref{eq:nuwp}) as a parametric LP problem, drawing on 
MDP theory, and to elucidate
the structure of its simplex tableaux.

Introducing \emph{variables} $x_j^a$ corresponding to occupancy measures 
$x_{ij}^{a, \tau}$ in Section \ref{s:sospgi}, we use
(\ref{eq:les}) to reformulate  (\ref{eq:nuwp})
as the following parametric LP problem:
\begin{equation}
\label{eq:nuwplp}
\begin{split}
& \max \,  \mathbf{x}^{1}
(\mathbf{R} - \nu \mathbf{1}) \\
& \text{subject to} \\
& \begin{bmatrix} \mathbf{x}^{0} &  \mathbf{x}^{1} \end{bmatrix}
\begin{bmatrix} 
(1-\beta) \mathbf{I} \\
\mathbf{I} - \beta \mathbf{P} 
\end{bmatrix} = \mathbf{e}_i \\
& \begin{bmatrix} \mathbf{x}^{0} &  \mathbf{x}^{1} \end{bmatrix}
 \geq \mathbf{0}.
\end{split}
\end{equation}
Notice that in LP (\ref{eq:nuwplp}) we write 
variables as a row  instead of as a column vector, which is
counter to
conventional usage in LP theory.
We do so for notational convenience, as in this way we avoid below using tranposed matrices such as $\transp{\mathbf{P}}$.

To analyze such an LP problem, we start by noticing that
 its \emph{basic feasible solutions}  (BFS) correspond to continuation sets 
$S \subseteq N$. We will thus refer to the \emph{$S$-active BFS}.
For each such $S$, we decompose 
the above vectors and matrices as 
\[
\mathbf{x}^{a} = \begin{bmatrix} \mathbf{x}_{S}^{a} &
  \mathbf{x}_{S^c}^{a} \end{bmatrix}, \quad 
\mathbf{e}_i = \begin{bmatrix} \mathbf{e}_{i S} &  \mathbf{e}_{i S^c}
  \end{bmatrix}, \quad 
\mathbf{P} = \begin{bmatrix} \mathbf{P}_{SS} & \mathbf{P}_{S S^c} \\
  \mathbf{P}_{S^c S} & \mathbf{P}_{S^c S^c} \end{bmatrix}, \quad
\mathbf{I} = \begin{bmatrix} \mathbf{I}_{S} & \mathbf{0}_{S S^c} \\
  \mathbf{0}_{S^c S} & \mathbf{I}_{S^c} \end{bmatrix},
\]
where, e.g., $\mathbf{I}_{S}$ is the identity matrix indexed by
$S \times S$, 
and  introduce the matrices 
\begin{equation}
\label{eq:matdef}
\begin{split}
\mathbf{P}^S & \triangleq \begin{bmatrix} \mathbf{P}_{S S} &
  \mathbf{P}_{S S^c} \\ \mathbf{0}_{S^c S} & \mathbf{I}_{S^c}
  \end{bmatrix}, \quad
\mathbf{P}^{S^c} \triangleq \begin{bmatrix} \mathbf{I}_{S} &
  \mathbf{0}_{S S^c} \\ \mathbf{P}_{S^c S} & \mathbf{P}_{S^c S^c}
  \end{bmatrix}, \\
\mathbf{B}^S & \triangleq \mathbf{I} - \beta \mathbf{P}^S, \quad
\mathbf{N}^S \triangleq \mathbf{I} - \beta \mathbf{P}^{S^c}, \quad
\mathbf{H}^S \triangleq \left\{\mathbf{B}^S\right\}^{-1}, \quad
\mathbf{A}^S \triangleq \mathbf{N}^S \mathbf{H}^S.
\end{split}
\end{equation}
Notice that $\mathbf{P}^S$ (resp. $\mathbf{P}^{S^c}$) is the transition-probability matrix
under the $S$-active (resp. $S^c$-active) policy.
Further, $\mathbf{B}^S$ is the \emph{basis matrix} in
(\ref{eq:nuwplp})
corresponding to 
the $S$-active BFS, whose \emph{basic variables} (corresponding to matrix
rows) are
$\begin{bmatrix} \mathbf{x}_S^1 & \mathbf{x}_{S^c}^0 \end{bmatrix}$ 
and $\mathbf{N}^S$ is the matrix of non-basic rows in
(\ref{eq:nuwplp}),
with associated \emph{non-basic variables}
$\begin{bmatrix} \mathbf{x}_S^0 & \mathbf{x}_{S^c}^1 \end{bmatrix}$.

We next use such a framework to represent performance
measures of concern under the $S$-active policy.
We first reformulate the constraints in LP (\ref{eq:nuwplp}) as

\begin{equation}
\label{eq:lesblocks}
\begin{bmatrix}
\begin{bmatrix} \mathbf{x}_{S}^{1} & \mathbf{x}_{S^c}^{0}
\end{bmatrix} &
\begin{bmatrix} \mathbf{x}_{S}^{0} & \mathbf{x}_{S^c}^{1}
\end{bmatrix}
\end{bmatrix}
\begin{bmatrix}
\mathbf{B}^S \\ \mathbf{N}^S
\end{bmatrix} = \mathbf{e}_i.
\end{equation}

We obtain the occupancies $x_{ij}^{a, S}$ under the $S$-active policy
as the corresponding 
$S$-active BFS, by setting to zero the non-basic variables, i.e., 
$\begin{bmatrix} \mathbf{x}_{iS}^{0, S} & \mathbf{x}_{i S^c}^{1, S}
\end{bmatrix} = \mathbf{0}$, and calculating the basic variables by 
solving the linear-equation system
\begin{equation}
\label{eq:xscalc}
\begin{bmatrix} \mathbf{x}_{S}^{1} & \mathbf{x}_{S^c}^{0}
\end{bmatrix} \mathbf{B}^S = \mathbf{e}_i \Longrightarrow 
\begin{bmatrix} \mathbf{x}_{i S}^{1, S} & \mathbf{x}_{i S^c}^{0, S}
\end{bmatrix} = \mathbf{e}_i \mathbf{H}^S.
\end{equation}

We now use (\ref{eq:xscalc}) to represent work measure (cf.\
 (\ref{eq:figilf})) vector $\mathbf{g}^S
= (g_i^S)_{i \in N}$:
\begin{equation}
\label{eq:gscalc}
g_i^S = \sum_{j \in S} x_{ij}^{1, S} = \begin{bmatrix} \mathbf{x}_{i S}^{1, S} & \mathbf{x}_{i S^c}^{0, S}
\end{bmatrix} 
\begin{bmatrix}
\mathbf{1}_S \\ \mathbf{0}_{S^c}
\end{bmatrix} 
= \mathbf{e}_i \mathbf{H}^S \begin{bmatrix}
\mathbf{1}_S \\ \mathbf{0}_{S^c}
\end{bmatrix}
 \Longrightarrow \mathbf{g}^S = \mathbf{H}^S \begin{bmatrix}
\mathbf{1}_S \\ \mathbf{0}_{S^c}
\end{bmatrix}.
\end{equation}

Similarly, we represent reward measure (cf.\
(\ref{eq:figilf})) vector $\mathbf{f}^S
= (f_i^S)_{i \in S}$:
\begin{equation}
\label{eq:fscalc}
f_i^S = \sum_{j \in S} R_j x_{ij}^{1, S}
= 
 \begin{bmatrix} \mathbf{x}_{i S}^{1, S} & \mathbf{x}_{i S^c}^{0, S}
\end{bmatrix} 
\begin{bmatrix}
\mathbf{R}_S \\ \mathbf{0}_{S^c}
\end{bmatrix} 
= \mathbf{e}_i \mathbf{H}^S \begin{bmatrix}
\mathbf{R}_S \\ \mathbf{0}_{S^c}
\end{bmatrix}  \Longrightarrow 
\mathbf{f}^S = \mathbf{H}^S \begin{bmatrix}
\mathbf{R}_S \\ \mathbf{0}_{S^c}
\end{bmatrix}.
\end{equation}

Further, we  represent
marginal workloads $\mathbf{w}^S =
(w_i^S)_{i \in N}$ using (\ref{eq:mwm}) and $\mathbf{g}_{S^c}^S = \mathbf{0}$:
\begin{equation}
\begin{split}
\label{eq:wsc1}
\mathbf{w}_S^S & = \mathbf{g}_S^S - \beta \mathbf{g}_S^S = (1-\beta)
\mathbf{g}_S^S \\
\mathbf{w}_{S^c}^S & = \mathbf{1}_{S^c} + \beta \mathbf{P}_{S^c N}
\mathbf{g}^S -  \mathbf{g}_{S^c}^S = \mathbf{1}_{S^c} + \beta \mathbf{P}_{S^c S}
\mathbf{g}_S^S.
\end{split}
\end{equation} 
Notice that (\ref{eq:wsc1}) implies that $\mathbf{w}^S > \mathbf{0}$,
i.e., \emph{marginal workloads are positive}, 
as mentioned before.
We now reformulate the first identities in (\ref{eq:wsc1}), using (\ref{eq:gscalc}), as 
\begin{equation}
\label{eq:wsc2}
\begin{bmatrix}
\mathbf{w}_S^S \\
- \mathbf{w}_{S^c}^S
\end{bmatrix}
= \mathbf{N}^S \mathbf{g}^S - 
\begin{bmatrix} \mathbf{0}_S \\ \mathbf{1}_{S^c} \end{bmatrix}
= \mathbf{N}^S \mathbf{H}^S \begin{bmatrix}
\mathbf{1}_S \\ \mathbf{0}_{S^c}
\end{bmatrix} - \begin{bmatrix} \mathbf{0}_S \\ \mathbf{1}_{S^c}
\end{bmatrix}
= \mathbf{A}^S  \begin{bmatrix}
\mathbf{1}_S \\ \mathbf{0}_{S^c}
\end{bmatrix} - \begin{bmatrix} \mathbf{0}_S \\ \mathbf{1}_{S^c}
\end{bmatrix}.
\end{equation}

Finally, we  obtain 
marginal reward  vector $\mathbf{r}^S =
(r_i^S)_{i \in N}$ from its definition in (\ref{eq:mrm}):
\begin{equation}
\begin{split}
\label{eq:rsc1}
\mathbf{r}_S^S & = \mathbf{f}_S^S  - \beta \mathbf{f}_S^S = (1-\beta) \mathbf{f}_S^S \\
\mathbf{r}_{S^c}^S & = \mathbf{R}_{S^c} + \beta \mathbf{P}_{S^c N}
\mathbf{f}^S -  \mathbf{f}_{S^c}^S = \mathbf{R}_{S^c} + \beta \mathbf{P}_{S^c N}
\mathbf{f}^S.
\end{split}
\end{equation}
We now reformulate (\ref{eq:rsc1}), using (\ref{eq:fscalc}), as 
\begin{equation}
\label{eq:rsc2}
\begin{bmatrix}
\mathbf{r}_S^S \\
- \mathbf{r}_{S^c}^S
\end{bmatrix}
= \mathbf{N}^S \mathbf{f}^S - 
\begin{bmatrix}  \mathbf{0}_S \\ \mathbf{R}_{S^c} \end{bmatrix}
= \mathbf{N}^S \mathbf{H}^S \begin{bmatrix}
\mathbf{R}_S \\ \mathbf{0}_{S^c}
\end{bmatrix} - 
\begin{bmatrix}  \mathbf{0}_S \\ \mathbf{R}_{S^c} \end{bmatrix}
= 
\mathbf{A}^S \begin{bmatrix}
\mathbf{R}_S \\ \mathbf{0}_{S^c}
\end{bmatrix} - 
\begin{bmatrix}  \mathbf{0}_S \\ \mathbf{R}_{S^c} \end{bmatrix}.
\end{equation}

Our next result characterizes the marginal workload and marginal
reward measures introduced above as reduced costs of corresponding
LP problems. It further characterizes the reduced costs of
parametric LP problem \textup{(\ref{eq:nuwplp})}.

\begin{proposition}
\label{pro:rclp}
Reduced costs for non-basic variables in the 
$S$-active BFS for LPs 
\[
\max \, \left\{\mathbf{x}^1 \mathbf{R}\colon 
\begin{bmatrix} \mathbf{x}^{0} &  \mathbf{x}^{1} \end{bmatrix}
\begin{bmatrix} 
(1-\beta) \mathbf{I} \\
\mathbf{I} - \beta \mathbf{P}
\end{bmatrix} = \mathbf{e}_i, \quad 
\begin{bmatrix} \mathbf{x}^{0} &  \mathbf{x}^{1} \end{bmatrix}
 \geq \mathbf{0}\right\},
\]
\[
\max \, \left\{\mathbf{x}^1 \mathbf{1}\colon 
\begin{bmatrix} \mathbf{x}^{0} &  \mathbf{x}^{1} \end{bmatrix}
\begin{bmatrix} 
(1-\beta) \mathbf{I} \\
\mathbf{I} - \beta \mathbf{P} 
\end{bmatrix} = \mathbf{e}_i, \quad 
\begin{bmatrix} \mathbf{x}^{0} &  \mathbf{x}^{1} \end{bmatrix}
 \geq \mathbf{0}\right\}
\]
and
\[
\max \, \left\{\mathbf{x}^1 (\mathbf{R} - \nu \mathbf{1})\colon 
\begin{bmatrix} \mathbf{x}^{0} &  \mathbf{x}^{1} \end{bmatrix}
\begin{bmatrix} 
(1-\beta) \mathbf{I} \\
\mathbf{I} - \beta \mathbf{P}
\end{bmatrix} = \mathbf{e}_i, \quad 
\begin{bmatrix} \mathbf{x}^{0} &  \mathbf{x}^{1} \end{bmatrix}
 \geq \mathbf{0}\right\}
\]
are given, respectively,  by \textup{(\ref{eq:rsc2})}, \textup{(\ref{eq:wsc2})},
 and 
\begin{equation}
\label{eq:rcvplp}
\begin{bmatrix}
\mathbf{r}_S^S - \nu \mathbf{w}_S^S \\
- \mathbf{r}_{S^c}^S + \nu \mathbf{w}_{S^c}^S
\end{bmatrix}.
\end{equation}
Therefore, such LPs' objectives can be
represented, respectively, as
\begin{equation}
\label{eq:rdl}
\mathbf{x}^1 \mathbf{R} = f_i^S - 
\sum_{j \in S} r_j^S x_j^0 +
 \sum_{j \in S^c} r_j^S x_j^1,
\end{equation}
\begin{equation}
\label{eq:wdl}
\mathbf{x}^1 \mathbf{1} = g_i^S - \sum_{j \in S} w_j^S x_j^0 +
 \sum_{j \in S^c} w_j^S x_j^1
\end{equation}
and
\begin{equation}
\label{eq:rwdl}
\mathbf{x}^1 (\mathbf{R} - \nu
\mathbf{1})
= f_i^S - \nu g_i^S - 
\sum_{j \in S} (r_j^S - \nu w_j^S) x_j^0 +
 \sum_{j \in S^c} (r_j^S - \nu w_j^S) x_j^1.
\end{equation}
\end{proposition}
\begin{proof}
The results follow 
from the representation of reduced costs
in LP theory, as given by (\ref{eq:wsc2}) and (\ref{eq:rsc2}),
along with the representation of an LP's objective 
in terms of the current BFS value and reduced costs.
Notice that the result for the LP with objective $\mathbf{x}^1
\mathbf{R}$ implies the results for the other LPs by appropriate choice
 of $\mathbf{R}$.
\end{proof}

The next result, which follows directly from Proposition
\ref{pro:rclp}, gives  representations of measures
$g_i^\tau$, $f_i^\tau$, and objective $f_i^{\tau} - \nu g_i^\tau$ 
relative to the $S$-active policy. Such
identities were first obtained in \citet{nmaap01} through algebraic arguments.

\begin{proposition}
\label{pro:dlaws} 
\mbox{ }
\begin{itemize}
\item[\textup{(a)}] $\displaystyle g_i^\tau = g_i^S - \sum_{j \in S} w_j^S x_{ij}^{0, \tau} +
 \sum_{j \in S^c} w_j^S x_{ij}^{1, \tau}$.
\item[\textup{(b)}] $\displaystyle f_i^\tau = f_i^S - \sum_{j \in S} r_j^S x_{ij}^{0,
 \tau} +
 \sum_{j \in S^c} r_j^S x_{ij}^{1, \tau}$.
\item[\textup{(c)}] $\displaystyle f_i^{\tau} - \nu g_i^\tau 
= f_i^S - \nu g_i^S - \sum_{j \in S} (r_j^S - \nu w_j^S) x_{ij}^{0, \tau} +
 \sum_{j \in S^c} (r_j^S - \nu w_j^S) x_{ij}^{1, \tau}$.
\end{itemize}
\end{proposition}

We next use the characterization of reduced costs in 
Proposition \ref{pro:rclp} to give a necessary and sufficient optimality test for 
the $S$-active BFS in parametric LP problem
\textup{(\ref{eq:nuwplp})}, and therefore for the $S$-active policy 
in $\nu$-wage problem 
(\ref{eq:nuwp}).

\begin{proposition}
\label{pro:lpoptest}
The $S$-active BFS is optimal for LP problem
\textup{(\ref{eq:nuwplp})} --- and hence
so is the $S$-active policy for $\nu$-wage problem 
\textup{(\ref{eq:nuwp})} --- for every initial state $i \in N$
iff 
\begin{equation}
\label{eq:nunujsineq}
\max \left\{\nu_j^S\colon j \in S^c\right\} 
\leq \nu \leq \min \left\{\nu_j^S\colon j \in S\right\}.
\end{equation}
\end{proposition}
\begin{proof}
The ``if'' part is the sufficient optimality test in LP theory
that checks nonnegativity of reduced costs for non-basic variables.
 The inequalities 
(\ref{eq:nunujsineq}) follow by reformulating such a condition, using
Proposition \ref{pro:rclp}(c),
positivity of marginal workloads $w_j^S$, and the
definition of marginal productivity rates $\nu_j^S$.

The ``only if'' part follows by considering a variation on LP
(\ref{eq:nuwplp})
where the right-hand side $\mathbf{e}_i$ is replaced by a positive
initial-state probability row vector $\mathbf{p} = (p_j)_{j \in N} >
\mathbf{0}$. The reduced-cost optimality test for such an LP is
the same as for the LPs with right-hand sides $\mathbf{e}_i$. Yet
such an LP is  non-degenerate, and hence 
the optimality test is also necessary for it. This completes the proof.
\end{proof}

\begin{table}[ht]
\caption{Augmented Tableau for $S$-active
  BFS, Ready for Pivoting on $a_{jj}^S$.}
\begin{center}
\begin{tabular}{cccccc} 
& $\mathbf{x}_S^1$ & $x_j^0$ & 
$\mathbf{x}_{S^c \setminus \{j\}}^0$  &  & 
\\  \cline{2-6} 
$\transp{\{\mathbf{x}_S^0\}}$ & 
\multicolumn{1}{|c}{$\mathbf{A}_{SS}^S$} &
  $\mathbf{A}_{Sj}^S$
 & \multicolumn{1}{c|}{$\mathbf{A}_{S, S^c \setminus \{j\}}^S$} & $\mathbf{w}_S^S$ & 
\multicolumn{1}{c|}{$\mathbf{r}_S^S$} \\
$x_j^1$ & 
\multicolumn{1}{|c}{$\mathbf{A}_{jS}^S$} & $\boxed{a_{jj}^S}$ & \multicolumn{1}{c|}{$\mathbf{A}_{j, S^c \setminus \{j\}}$} & 
$-w_j^S$ & \multicolumn{1}{c|}{$-r_j^S$} \\
$\transp{\{\mathbf{x}_{S^c \setminus \{j\}}^1\}}$ & 
\multicolumn{1}{|c}{$\mathbf{A}_{S^c \setminus \{j\}, S}^S$} & 
$\mathbf{A}_{S^c \setminus \{j\}, j}^S$ & 
\multicolumn{1}{c|}{$\mathbf{A}_{S^c \setminus \{j\}, S^c \setminus \{j\}}^S$} & 
$- \mathbf{w}_{S^c \setminus \{j\}}^S$ &   
\multicolumn{1}{c|}{$- \mathbf{r}_{S^c \setminus \{j\}}^S$} \\
$x_j^0$ & \multicolumn{1}{|c}{$\mathbf{0}_{jS}$} & $1$ & \multicolumn{1}{c|}{$\mathbf{0}_{j, S^c \setminus
  \{j\}}$} & $0$ & \multicolumn{1}{c|}{$0$} \\ 
\cline{2-6} 
\end{tabular}
\end{center}
\label{tab:invtabl}
\end{table}

We now have all the elements to represent the 
\emph{parametric simplex tableau}  for the 
$S$-active BFS, as shown in Table \ref{tab:invtabl}.
Notice that such tableau  is 
\emph{transposed} relative to conventional simplex tableaux, so that
non-basic variables $\mathbf{x}_S^0$ and $\mathbf{x}_{S^c}^1$
correspond to \emph{rows}, whereas basic variables $\mathbf{x}_S^1$ and 
$\mathbf{x}_{S^c}^0$ correspond to \emph{columns}.
Further, it includes two columns of reduced costs for 
non-basic variables,  corresponding to the first two LP problems in 
Proposition \ref{pro:rclp}.

Actually, Table \ref{tab:invtabl} represents an \emph{augmented
  tableau}, ready for \emph{pivoting} on element $a_{jj}^S$, with 
$j \in S^c$, i.e.,
for taking variable $x_j^0$ out of the current basis, and 
putting  $x_j^1$ into the basis.
To prepare the ground for carrying out
 the required pivoting operations, the tableau includes
an extra row, corresponding to basic variable $x_j^0$.

After such a pivoting step is performed, one obtains
the tableau for the $S \cup \{j\}$-active
  BFS, shown in Tables \ref{tab:tabl2} and \ref{tab:tabl3}.

\begin{table}[ht]
\caption{Part $\mathbf{A}^{S \cup \{j\}}$ of 
Tableau for $S \cup \{j\}$-Active
  BFS, Obtained by Pivoting.}
\begin{center}
\begin{tabular}{cccc} 
& $\mathbf{x}_S^1$ & $x_j^1$ & 
$\mathbf{x}_{S^c \setminus \{j\}}^0$   
\\  \cline{2-2} \cline{3-3} \cline{4-4}
$\transp{\{\mathbf{x}_S^0\}}$ & 
\multicolumn{1}{|c}{$\mathbf{A}_{SS}^S - \frac{\mathbf{A}_{Sj}^S
    \mathbf{A}_{jS}^S}{a_{jj}^S}$} &
  $\frac{\mathbf{A}_{Sj}^S}{a_{jj}^S}$
 & \multicolumn{1}{c|}{$\mathbf{A}_{S, S^c \setminus \{j\}}^S - 
   \frac{\mathbf{A}_{Sj}^S \mathbf{A}_{j, S^c \setminus \{j\}}^S}{a_{jj}^S}$}  \\
$x_j^0$ & 
\multicolumn{1}{|c}{$-\frac{\mathbf{A}_{jS}^S}{a_{jj}^S}$} & $\frac{1}{a_{jj}^S}$ & \multicolumn{1}{c|}{$-\frac{\mathbf{A}_{j, S^c \setminus \{j\}}^S}{a_{jj}^S}$} \\
$\transp{\{\mathbf{x}_{S^c \setminus \{j\}}^1\}}$ & 
\multicolumn{1}{|c}{$\mathbf{A}_{S^c \setminus \{j\}, S}^S - 
 \frac{\mathbf{A}_{S^c \setminus \{j\}, j}^S
   \mathbf{A}_{jS}^S}{a_{jj}^S}$} & 
$\frac{\mathbf{A}_{S^c \setminus \{j\}, j}^S}{a_{jj}^S}$ & 
\multicolumn{1}{c|}{$\mathbf{A}_{S^c \setminus \{j\}, S^c \setminus \{j\}}^S - 
\frac{\mathbf{A}_{S^c \setminus \{j\}, j}^S \mathbf{A}_{j, S^c \setminus \{j\}}^S}{a_{jj}^S}$}  \\
\cline{2-2} \cline{3-3} \cline{4-4}
\end{tabular}
\end{center}
\label{tab:tabl2}
\end{table}

\begin{table}[ht]
\caption{Parts $\mathbf{w}^{S \cup \{j\}}$, 
$\mathbf{r}^{S \cup \{j\}}$ of 
Tableau for $S \cup \{j\}$-Active
  BFS, Obtained by Pivoting.}
\begin{center}
\begin{tabular}{ccc} 
\cline{2-2} \cline{3-3} 
$\transp{\{\mathbf{x}_S^0\}}$ & \multicolumn{1}{|c}{$\mathbf{w}_S^S + 
 \frac{w_j^S}{a_{jj}^S} \mathbf{A}_{Sj}^S$} & 
\multicolumn{1}{c|}{$\mathbf{r}_S^S + 
 \frac{r_j^S}{a_{jj}^S}  \mathbf{A}_{Sj}^S$} \\
$x_j^0$ & \multicolumn{1}{|c}{$\frac{w_j^S}{a_{jj}^S}$} & \multicolumn{1}{c|}{$\frac{r_j^S}{a_{jj}^S}$} \\
$\transp{\{\mathbf{x}_{S^c \setminus \{j\}}^1\}}$ & 
\multicolumn{1}{|c}{$- \mathbf{w}_{S^c \setminus \{j\}}^S + 
 \frac{w_j^S}{a_{jj}^S} \mathbf{A}_{S^c \setminus \{j\}, j}^S$} &   
\multicolumn{1}{c|}{$- \mathbf{r}_{S^c \setminus \{j\}}^S + 
\frac{r_j^S}{a_{jj}^S}  \mathbf{A}_{S^c \setminus \{j\}, j}^S$} \\
\cline{2-2} \cline{3-3}
\end{tabular}
\end{center}
\label{tab:tabl3}
\end{table}

We must further  elucidate the structure of the
initial tableau, corresponding to the 
$\emptyset$-active BFS. 
Letting $S = \emptyset$, we readily obtain from 
(\ref{eq:matdef}) that
\begin{equation}
\label{eq:mats0}
\mathbf{B}^\emptyset = (1-\beta) \mathbf{I}, \quad
\mathbf{N}^\emptyset = \mathbf{I} - \beta \mathbf{P}, \quad
\mathbf{H}^\emptyset = \frac{1}{1-\beta} \mathbf{I}, \quad
\mathbf{A}^\emptyset = \frac{1}{1-\beta} (\mathbf{I} - \beta \mathbf{P}).
\end{equation}
Further, using (\ref{eq:mats0}), (\ref{eq:wsc2}), and
(\ref{eq:rsc2}), we obtain the initial reduced costs from
\begin{equation}
\label{eq:initrc}
\mathbf{w}^\emptyset = \mathbf{1} \quad \text{and} \quad
\mathbf{r}^\emptyset = \mathbf{R}.
\end{equation}

We are now ready to formulate the \emph{conventional-pivoting} (CP) index 
algorithm, shown in Table \ref{tab:cpgia}, which 
implements  the parametric
simplex approach of \citet{kallenb86}.  
In Table \ref{tab:cpgia} we have adopted an algorithm-like
notation, substituting superscript counters $(k)$ for superscript sets $S_k$.
Notice that such an algorithm only applies to the  discounted
case $\beta < 1$, since the computation of $\mathbf{A}^{(0)}$ involves
division by $1-\beta$.

We next assess the computational complexity of the CP algorithm. 

\begin{proposition}
\label{pro:cccpa}
The CP algorithm performs
$2n^{3}+ O(n^2)$ arithmetic operations.
\end{proposition}
\begin{proof}
Counting shows that
the algorithm performs $n^{3}
+ O(n^2)$ multiplications and divisions and the same order of 
additions and subtractions, which yields the result.
\end{proof}

\begin{table}[tb]
\caption{The \textbf{CP} Gittins-Index Algorithm.}
\begin{center}
\fbox{%
\begin{minipage}{\textwidth}
\textbf{ALGORITHM CP:} 
\begin{tabbing}
$S_0 := \emptyset$,  
$\displaystyle \mathbf{A}^{(0)} := \frac{1}{1-\beta}(\mathbf{I} - \beta \mathbf{P})$, 
$\mathbf{w}^{(0)} :=
\mathbf{1}$,  $\displaystyle \mathbf{r}^{(0)} := 
 \mathbf{R}$
\end{tabbing}

\begin{tabbing}
\textbf{for} \= $k := 1$ \textbf{to} $n$ \textbf{do} \\
 \> \textbf{pick} 
 $i_{k} \in \argmax
      \left\{r^{({k-1})}_i/w^{({k-1})}_i\colon
                i \in S_{k-1}^c\right\}; \,\, \nu_{i_k}^* := 
 r^{({k-1})}_{i_k}/w^{({k-1})}_{i_k}, \,\, S_{k} := S_{k-1} \cup \{i_{k}\}$  \\
 \> \textbf{if} \= $k < n$ \textbf{then} \\
 \>  \> $p^{(k)} := a_{i_k i_k}^{(k-1)}, \,\, a_{i_k i_k}^{(k-1)} := 1, \,\, 
\mathbf{v}^{(k)} := (1/p^{(k)}) \mathbf{A}_{N i_k}^{(k-1)}, \,\, 
\mathbf{h}^{(k)} := -\mathbf{A}_{i_k N}^{(k-1)}$
 \\ \\
 \>  \> $\mathbf{A}^{(k)} := \mathbf{A}^{(k-1)} + \mathbf{v}^{(k)}
 \mathbf{h}^{(k)}, \,\, \mathbf{A}_{N i_k}^{(k)} := \mathbf{v}^{(k)}, \,\,
\mathbf{A}_{i_k N}^{(k)} := (1/p^{(k)}) \mathbf{h}^{(k)}$
\\ \\
 \>   \>  $\mathbf{w}_{S_k^c}^{(k)} := \mathbf{w}_{S_k^c}^{({k-1})} - 
  w_{i_k}^{({k-1})} \mathbf{A}_{S_k^c i_k}^{({k})}, \,\,
\mathbf{w}_{S_{k-1}}^{(k)} := \mathbf{w}_{S_{k-1}}^{({k-1})} + 
  w_{i_k}^{({k-1})} \mathbf{A}_{S_{k-1} i_k}^{({k})}, \,\, 
w_{i_k}^{(k)} := w_{i_k}^{(k-1)}/p^{(k)}$ \\ \\
\>  \>   $\mathbf{r}_{S_k^c}^{(k)} := \mathbf{r}_{S_k^c}^{({k-1})} - 
  r_{i_k}^{({k-1})} \mathbf{A}_{S_k^c i_k}^{({k})}, \,\,
\mathbf{r}_{S_{k-1}}^{(k)} := \mathbf{r}_{S_{k-1}}^{({k-1})} + 
  r_{i_k}^{({k-1})} \mathbf{A}_{S_{k-1} i_k}^{({k})}, \,\, 
r_{i_k}^{(k)} := r_{i_k}^{(k-1)}/p^{(k)}$ \\
 \>  \textbf{end } \{ if \} \\
\textbf{end} \{ for \}
\end{tabbing}
\end{minipage}}
\end{center}
\label{tab:cpgia}
\end{table}

\section{Exploiting Special Structure}
\label{s:ess}
We set out in this section to exploit special structure to 
reduce the number of operations performed in a pivoting
step.

\subsection{Updating Marginal Productivity Rates}
\label{s:umpr}
We start by showing that there is no need to update marginal 
rewards $r_j^S$ in the tableaux. It suffices to update 
required marginal work measures $w_j^S$, and then use them
to update
required marginal productivity rates $\nu_j^S$.
While the next result is proven in 
 \citet{nmmp02}, we
present here for self-completeness a new proof, drawing on the pivot step 
that leads from the tableau in Table \ref{tab:invtabl} to that in 
Tables \ref{tab:tabl2}--\ref{tab:tabl3}.

\begin{proposition}
\label{pro:mprupd} For $i \in N$ and $j \in S^c$,
$\displaystyle \nu_i^{S \cup \{j\}} = \nu_j^S - \frac{w_i^S}{w_i^{S \cup \{j\}}}
(\nu_j^S - \nu_i^S)$.
\end{proposition}
\begin{proof}
Start with the case $i \in S$. 
From Table \ref{tab:tabl3} and $\nu_i^{S \cup \{j\}} = 
r_i^{S \cup \{j\}}/w_i^{S \cup \{j\}}$, we have
\begin{align*}
\nu_i^{S \cup \{j\}} & = 
\frac{\nu_i^S w_i^S + \nu_j^S \frac{w_j^S}{a_{jj}^S} a_{ij}^S}{w_i^S +
  \frac{w_j^S}{a_{jj}^S} a_{ij}^S} 
= \nu_j^S - \frac{w_i^S}{w_i^S + \frac{w_j^S}{a_{jj}^S} a_{ij}^S} 
(\nu_j^S - \nu_i^S) = 
\nu_j^S - \frac{w_i^S}{w_i^{S \cup \{j\}}}
(\nu_j^S - \nu_i^S).
\end{align*}

Consider now the case $i = j$. Then, again Table \ref{tab:tabl3} yields
 that
\begin{align*}
\nu_j^{S \cup \{j\}} = 
\frac{r_j^S/a_{jj}^S}{w_j^S/a_{jj}^S} = \frac{r_j^S}{w_j^S} = \nu_j^{S}.
\end{align*}

Finally, consider the case $i \in S^c \setminus \{j\}$. Then, using
again Table \ref{tab:tabl3}
we obtain that
\begin{align*}
\nu_i^{S \cup \{j\}} & =
\frac{\nu_i^S w_i^S - \nu_j^S \frac{w_j^S}{a_{jj}^S} a_{ij}^S}{w_i^S -
  \frac{w_j^S}{a_{jj}^S} a_{ij}^S} 
= \nu_j^S - \frac{w_i^S}{w_i^S - \frac{w_j^S}{a_{jj}^S} a_{ij}^S} 
(\nu_j^S - \nu_i^S) = 
\nu_j^S - \frac{w_i^S}{w_i^{S \cup \{j\}}}
(\nu_j^S - \nu_i^S),
\end{align*}
as required. This completes the proof.
\end{proof}

\subsection{The Reduced Tableau}
\label{s:rpmas}
We proceed to elucidate which is the minimal 
information that needs be updated at each pivoting step, 
which will be stored in a \emph{reduced tableau}.
For such a purpose, we first observe that 
basis inverse matrix $\mathbf{H}^S$ is readily partitioned and
represented as
\begin{equation}
\label{eq:hsdec}
\begin{split}
\mathbf{H}^S & = 
\begin{bmatrix}
\mathbf{H}_{SS}^S & \mathbf{H}_{SS^c}^S \\
\mathbf{H}_{S^c S}^S & \mathbf{H}_{S^c S^c}^S
\end{bmatrix} = 
\begin{bmatrix}
(\mathbf{I}_{S} - \beta \mathbf{P}_{SS})^{-1}
  & \frac{\beta}{1-\beta} (\mathbf{I}_{S} - \beta \mathbf{P}_{SS})^{-1}
\mathbf{P}_{S S^c} \\
\mathbf{0}_{S^c S} & \frac{1}{1-\beta} \mathbf{I}_{S^c}
\end{bmatrix} \\
& = 
\begin{bmatrix}
\mathbf{H}_{SS}^S  & \frac{\beta}{1-\beta} \mathbf{H}_{SS}^S
\mathbf{P}_{S S^c} \\
\mathbf{0}_{S^c S} & \frac{1}{1-\beta} \mathbf{I}_{S^c}
\end{bmatrix}.
\end{split}
\end{equation}

We now use (\ref{eq:hsdec}) to partition and represent matrix
$\mathbf{A}^S$ as
\begin{equation}
\label{eq:aspart}
\begin{split}
\mathbf{A}^S & = \begin{bmatrix}
\mathbf{A}_{SS}^S & \mathbf{A}_{SS^c}^S \\
\mathbf{A}_{S^c S}^S & \mathbf{A}_{S^c S^c}^S
\end{bmatrix} = \mathbf{N}^S \mathbf{H}^S = \begin{bmatrix}
(1-\beta) \mathbf{I}_{S} & \mathbf{0}_{SS^c} \\
-\beta \mathbf{P}_{S^cS} & 
\mathbf{I}_{S^c} - \beta \mathbf{P}_{S^c S^c}
\end{bmatrix} 
\begin{bmatrix}
\mathbf{H}_{SS}^S  & \frac{\beta}{1-\beta} \mathbf{H}_{SS}^S
\mathbf{P}_{S S^c} \\
\mathbf{0}_{S^c S} & \frac{1}{1-\beta} \mathbf{I}_{S^c}
\end{bmatrix}
 \\
& = 
\begin{bmatrix}
(1-\beta) \mathbf{H}_{SS}^S & 
\beta \mathbf{H}_{SS}^S \mathbf{P}_{SS^c} \\
-\beta \mathbf{P}_{S^cS} \mathbf{H}_{SS}^S  & 
\frac{1}{1-\beta} 
(\mathbf{I}_{S^c} - \beta \mathbf{P}_{S^cS^c}) -\frac{\beta^2}{1-\beta} \mathbf{P}_{S^cS}
\mathbf{H}_{SS}^S
\mathbf{P}_{SS^c}
\end{bmatrix} \\
& = \begin{bmatrix}
(1-\beta) \mathbf{H}_{SS}^S & \mathbf{A}_{SS^c}^S \\
\mathbf{A}_{S^c S}^S & \frac{1}{1-\beta} 
(\mathbf{I}_{S^c} - \beta \mathbf{P}_{S^cS^c}
+\beta \mathbf{A}_{S^cS}^S \mathbf{P}_{SS^c})
\end{bmatrix}.
\end{split}
\end{equation}

From (\ref{eq:aspart}) and the above discussion, it should be clear
that it suffices to keep the
 information on the $S$-active BFS
shown in the reduced tableau in Table \ref{tab:redtabl}.

\begin{table}[ht]
\caption{Reduced Tableau for $S$-Active
  BFS.}
\begin{center}
\begin{tabular}{cccc} 
& $\mathbf{x}_S^1$ &    & 
\\  \cline{2-4} 
$\transp{\{\mathbf{x}_{S^c}^1\}}$ & 
\multicolumn{1}{|c|}{$\mathbf{A}_{S^cS}^S$}
 & 
$\mathbf{w}_{S^c}^S$ &   
\multicolumn{1}{c|}{$\boldsymbol{\nu}_{S^c}^S$}  \\
\cline{2-4} 
\end{tabular}
\end{center}
\label{tab:redtabl}
\end{table}

To obtain the reduced tableau for the $S \cup \{j\}$-active
  BFS from that in Table \ref{tab:redtabl}, we first 
compute pivot element $a_{jj}^S$ by noting that, by (\ref{eq:aspart}), 
\begin{equation}
\label{eq:ajjs}
a_{jj}^S = \frac{1 - \beta p_{jj} + \beta  \mathbf{A}_{jS}^S \mathbf{P}_{Sj}}{1-\beta}.
\end{equation}

We then use again (\ref{eq:aspart}) to compute the required vector
\begin{equation}
\label{eq:arvect}
\begin{split}
\mathbf{A}_{S^c \setminus \{j\}, j}^S & = - \frac{\beta}{1-\beta}
\left\{\mathbf{P}_{S^c \setminus \{j\}, j} - 
  \mathbf{A}_{S^c \setminus \{j\}, S}^S \mathbf{P}_{S j}\right\}.
\end{split}
\end{equation}

We can then readily compute the reduced tableau for 
the $S \cup \{j\}$-active BFS, as shown in Tables 
\ref{tab:rtscj1}--\ref{tab:rtscj2}.

\begin{table}[ht]
\caption{Reduced Tableau for $S \cup \{j\}$-Active
  BFS: $\mathbf{A}^{S \cup \{j\}}$ Part.}
\begin{center}
\begin{tabular}{ccc} 
& $\mathbf{x}_S^1$ & $x_j^1$ 
\\  \cline{2-3}
$\transp{\{\mathbf{x}_{S^c \setminus \{j\}}^1\}}$ & 
\multicolumn{1}{|c}{$\mathbf{A}_{S^c \setminus \{j\}, S}^S - 
 \frac{\mathbf{A}_{S^c \setminus \{j\}, j}^S
   \mathbf{A}_{jS}^S}{a_{jj}^S}$} & 
\multicolumn{1}{c|}{$\frac{\mathbf{A}_{S^c \setminus \{j\}, j}^S}{a_{jj}^S}$} 
   \\
\cline{2-3}
\end{tabular}
\end{center}
\label{tab:rtscj1}
\end{table}

\begin{table}[ht]
\caption{Reduced Tableau for $S \cup \{j\}$-Active
  BFS:  $\mathbf{w}^{S \cup \{j\}}$ and $\boldsymbol{\nu}^{S \cup
    \{j\}}$ Parts.}
\begin{center}
\begin{tabular}{ccc} 
\cline{2-2} \cline{3-3} 
$\transp{\{\mathbf{x}_{S^c \setminus \{j\}}^1\}}$ & 
\multicolumn{1}{|c}{$\mathbf{w}_{S^c \setminus \{j\}}^S -
 \frac{w_j^S}{a_{jj}^S} \mathbf{A}_{S^c \setminus \{j\}, j}^S$} &   
\multicolumn{1}{c|}{$\nu_j^S - \frac{w_i^S}{w_i^{S \cup \{j\}}}
(\nu_j^S - \nu_i^S), \, i \in S^c \setminus \{j\}$} \\
\cline{2-2} \cline{3-3}
\end{tabular}
\end{center}
\label{tab:rtscj2}
\end{table}

\section{The Fast-Pivoting Algorithm}
\label{s:piat}
Drawing on the above,  
we readily obtain the
\emph{fast-pivoting} algorithm FP(EO)  in Table \ref{tab:fasteo},
where again we write, e.g., $w_i^{S_k}$ as $w_i^{(k)}$.
The input $\mathrm{EO}$ is a Boolean variable than, when fed the
value \texttt{true}/1, makes the algorithm produce an \emph{extended
  output}.
Such a capability is useful
in some settings, where one may need, besides the 
Gittins index $\nu_i^*$, marginal work and marginal productivity 
measures $w_i^{S_k}$ and $\nu_i^{S_k}$ corresponding to the continuation sets 
$S_k$ generated by the algorithm. 
Thus,  e.g., such quantities are 
used in \citet{nmjoc06} as input for a fast algorithm that computes the \emph{switching
  index} of \citet{asatene96}.

We next assess the computational complexity of the FP(EO)
algorithm.

\begin{proposition}
\label{pro:ccfp} \text{ }
\begin{itemize}
\item[\textup{(a)}]
The \textup{FP(0)} algorithm 
performs $(2/3)n^{3}+ O(n^2)$ arithmetic operations.
\item[\textup{(b)}]
The \textup{FP(1)} algorithm 
performs $(4/3)n^{3}+ O(n^2)$ arithmetic operations.
\end{itemize}
\end{proposition}
\begin{proof}
(a) Step $k$ of the algorithm performs about $(k+1)(2(n-k)+1)$
multiplications and divisions, yielding a total of
\[
\sum_{k=2}^{n-1} \left( k+1\right) \left( 2(n-k)+1\right) + O(n)
=(1/3)n^{3}+(3/2)n^{2} + O(n).
\] 
The same count is obtained for additions and subtractions,
which yields the result.

(b) 
Counting shows that the algorithm performs 
$(2/3) n^{3}+ (5/2) n^{2} + O(n)$ multiplications and divisions, as well
as additions and subtractions, which gives the result.
\end{proof}

Instead of computing marginal productivity rates
$\nu_i^S$ as stated, we could have updated marginal rewards $r_i^S$, 
using the update formula 
\[
\mathbf{r}_{S_k^c}^{S_k} := \mathbf{r}_{S_k^c}^{S_{k-1}} - 
  r_{i_k}^{S_{k-1}} \mathbf{A}_{S_k^c i_k}^{S_{k}},
\]
and then have computed $\nu_i^{S_k} = r_i^{S_k}/w_i^{S_k}$ for $i \in
S_k^c$.
The complexity of such an alternative scheme is the same as that of the 
one proposed.

Notice that the FP(0) algorithm applies both to the
discounted and the undiscounted index. The latter is computed by
setting $\beta = 1$ in the stated calculations.

\begin{table}[tb]
\caption{The \textbf{FP}($\mathrm{EO}$) Gittins-Index Algorithm.}
\begin{center}
\fbox{%
\begin{minipage}{\textwidth}
\textbf{ALGORITHM FP}($\mathrm{EO}$): 
\begin{tabbing}
$S_0 := \emptyset$,  $\boldsymbol{w}^{(0)} :=
\mathbf{1}$, $\boldsymbol{\nu}^{(0)} := \mathbf{R}$
\end{tabbing}

\begin{tabbing}
\textbf{for} \= $k := 1$ \textbf{to} $n$ \textbf{do} \\
 \> \textbf{pick} 
 $\displaystyle i_{k} \in \argmax
      \left\{\nu^{({k-1})}_i\colon
                i \in S_{k-1}^c\right\}; \quad$
 $\displaystyle \nu_{i_k}^* := 
 \nu^{({k-1})}_{i_k}$, $S_{k} := S_{k-1} \cup \{i_{k}\}$  \\
 \> \textbf{if} \= $k = 1$ \textbf{then} \\
 \> \> $\alpha^{(1)} := -\beta/(1 - \beta p_{i_1
 i_1})$,
$\mathbf{A}_{S_1^c i_1}^{(1)} := 
  \alpha^{(1)} \mathbf{P}_{S_1^c i_1}$ \\
 \> \> \textbf{if } $\mathrm{EO}$ \textbf{then} $\mathbf{A}_{i_1 S_1^c}^{(1)} := -\alpha^{(1)} \mathbf{P}_{i_1
 S_1^c}$ \textbf{end } \{ if \} \\
 \> \textbf{else if}  $k < n$ \textbf{or} \{$\mathrm{EO}$ \textbf{and} $k = n$\} \textbf{then} \\
 \>  \>  $\alpha^{(k)} := -\beta/\{1 - \beta (p_{i_k
 i_k} - \mathbf{A}_{i_k S_{k-1}}^{({k-1})} \mathbf{P}_{S_{k-1}
 i_k})\}$ \\
\> \> 
$\mathbf{A}_{S_k^c i_k}^{(k)}
 := \alpha^{(k)} \{\mathbf{P}_{S_k^c i_k} -
 \mathbf{A}_{S_k^c S_{k-1}}^{(k-1)} \mathbf{P}_{S_{k-1}
 i_k}\}$, \, $\mathbf{A}_{S_k^c S_{k-1}}^{(k)}
 := \mathbf{A}_{S_k^c S_{k-1}}^{({k-1})} - \mathbf{A}_{S_{k}^c 
 i_k}^{({k})} \mathbf{A}_{i_k S_{k-1}}^{({k-1})}$
 \\
 \> \> \textbf{if } \= $\mathrm{EO}$ \textbf{then} \\
\>  \>  \>
$\mathbf{A}_{i_k S_k^c}^{(k)} = -\alpha^{(k)} \{\mathbf{P}_{i_k S_k^c} +
  \mathbf{P}_{i_k S_{k-1}} \mathbf{A}_{S_{k-1} S_k^c}^{(k-1)}\}$, \,
 $\mathbf{A}_{S_{k-1} S_k^c}^{(k)} := 
 \mathbf{A}_{S_{k-1} S_k^c}^{(k-1)} + \mathbf{A}_{S_{k-1} i_k}^{(k-1)}
 \mathbf{A}_{i_k S_k^c}^{(k)}$
\\
 \>  \>  \textbf{end } \{ if \} \\
 \>  \textbf{end } \{ if \} \\
 \>  $\mathbf{w}_{S_k^c}^{(k)} := 
\mathbf{w}_{S_k^c}^{(k-1)} - 
  w_{i_k}^{({k-1})} \mathbf{A}_{S_k^c i_k}^{({k})}$, $\displaystyle \nu_{i}^{(k)} := 
 \nu_{i_k}^* - \{w_i^{({k-1})}/w_i^{(k)}\} \{\nu_{i_k}^* -
 \nu_i^{({k-1})}\}, i \in S_k^c$ \\
 \> \textbf{if } \= $\mathrm{EO}$ \textbf{then} \\
 \>  \>
$w_{i_k}^{({k})} := -\{\alpha^{(k)}(1-\beta)/\beta\}
 w_{i_k}^{({k-1})}$, \, $\mathbf{w}_{S_{k-1}}^{(k)} := \mathbf{w}_{S_{k-1}}^{(k-1)} + 
 w_{i_k}^{({k})} \mathbf{A}_{S_{k-1} i_k}^{({k-1})}$
\\
 \> \> $\nu_{i}^{(k)} := 
 \nu_{i_k}^* - \{w_i^{({k-1})}/w_i^{(k)}\} \{\nu_{i_k}^* -
 \nu_i^{({k-1})}\}, i \in S_{k-1}$, \, $\nu_{i_k}^{(k)} :=
 \nu_{i_k}^*$ \\
 \>  \textbf{end } \{ if \} \\
\textbf{end} \{ for \}
\end{tabbing}
\end{minipage}}
\end{center}
\label{tab:fasteo}
\end{table}

\section{Relation with the State-Elimination Algorithm}
\label{s:rsea}
 \citet{kattaset04} and \citet{sonin05}
 introduced Gittins-index algorithms based on 
state-elimination ideas, discussing them
 as they apply to  bandits with
a terminal state. Yet, any bandit with a positive
discount factor is readily transformed into such case. See \citet{sonin05}.
Table \ref{tab:elim} formulates the \emph{state-elimination} (SE)
algorithm.

The above discussions allow us to show
that such algorithm computes  the Gittins index, by elucidating
its
relation with
 the CP algorithm.

\begin{table}[tb]
\caption{The \textbf{SE} Gittins-index algorithm.}
\begin{center}
\fbox{%
\begin{minipage}{\textwidth}
\textbf{ALGORITHM SE:} 
\begin{tabbing}
$S_0 := \emptyset, \,\, \tilde{\mathbf{P}}^{(0)} := \beta \mathbf{P},
\,\, \boldsymbol{\beta}^{{(0)}} :=
\beta \mathbf{1}, \,\, \tilde{\mathbf{r}}^{(0)} := (1-\beta) \mathbf{R}$
\end{tabbing}

\begin{tabbing}
\textbf{for} \= $k := 1$ \textbf{to} $n$ \textbf{do} \\
 \> \textbf{pick} 
 $i_{k} \in \argmax
      \left\{\tilde{r}^{({k-1})}_i/\{1-\beta_i^{(k-1)}\}\colon
                i \in S_{k-1}^c\right\}$ \\
 \>  $\nu_{i_k}^* := \tilde{r}^{({k-1})}_{i_k}/\{1-\beta_{i_k}^{(k-1)}\}, \, \, S_{k} := S_{k-1} \cup \{i_{k}\}$  \\
 \> \textbf{if} \= $k < n$ \textbf{then} \\
 \>  \>  $\tilde{\mathbf{P}}_{S_k^c i_k}^{(k)} :=
 \{1/(1-\tilde{p}_{i_k i_k}^{(k-1)})\} \tilde{\mathbf{P}}_{S_k^c
 i_k}^{(k-1)}, \,\, \tilde{\mathbf{P}}_{S_k^c S_{k}^c}^{(k)}
 := \tilde{\mathbf{P}}_{S_k^c S_{k}^c}^{({k-1})} +  \tilde{\mathbf{P}}_{S_{k}^c 
 i_k}^{({k})} \tilde{\mathbf{P}}_{i_k S_{k}^c}^{({k-1})}$
\\ \\
 \>  \>   $\boldsymbol{\beta}_{S_k^c}^{(k)} :=
 \tilde{\mathbf{P}}_{S_k^c S_k^c}^{(k)} \mathbf{1}_{S_k^c}, \,\,
\tilde{\mathbf{r}}_{S_k^c}^{(k)} :=
 \tilde{\mathbf{r}}_{S_k^c}^{(k-1)} + 
\tilde{r}_{i_k}^{(k-1)} \tilde{\mathbf{P}}_{S_k^c i_k}^{(k)}$ \\
 \>  \textbf{end } \{ if \} \\
\textbf{end} \{ for \}
\end{tabbing}
\end{minipage}}
\end{center}
\label{tab:elim}
\end{table}

\begin{proposition}
\label{pro:stelequiv}
The SE algorithm computes the Gittins index.
\end{proposition}
\begin{proof}
Straightforward algebra
 yields that the quantities computed in the
SE algorithm are related to those 
computed in the CP algorithm
 by 
\begin{equation}
\label{eq:pivrel}
\begin{split}
\tilde{\mathbf{P}}_{S_k^c S_{k}^c}^{(k)} & = \mathbf{I} - (1-\beta) \mathbf{A}_{S_k^c S_{k}^c}^{(k)}  \\
\tilde{r}_{S_k^c}^{(k)} & = (1-\beta) \mathbf{r}_{S_k^c}^{(k)} \\
\boldsymbol{\beta}_{S_k^c}^{(k)} & = \mathbf{1}_{S_k^c} - (1-\beta) \mathbf{w}_{S_k^c}^{(k)},
\end{split}
\end{equation}
for each $k$.
Therefore, the SE algorithm computes the same index as does the 
CP  algorithm. This completes the proof.
\end{proof}

We next assess the computational complexity of the 
SE algorithm.

\begin{proposition}
\label{pro:ccea}
The SE algorithm 
performs $n^{3}+ O(n^2)$ arithmetic operations.
\end{proposition}
\begin{proof}
Step $k$ of the algorithm performs $3(n-k)+(n-k)^{2} + O(1)$
 multiplications and divisions, giving a total of
\[
n^{2}+\sum_{k=1}^{n-1}\left\{3(n-k)+(n-k)^{2}\right\} + O(n)
= (1/3)n^{3}+2n^{2}+ O(n).
\] 

It further performs $2(n-k)^{2}+2(n-k) + O(1)$ additions and
subtractions, giving 
\[
2 \sum_{k=1}^{n-1}\left\{(n-k)^{2}+(n-k)\right\} + O(1)
=(2/3)n^{3}+ O(n).
\]
The stated total operation count follows.
\end{proof}

\section{Computation of the Undiscounted Gittins Index}
\label{s:cargi}
To motivate the interest of computing the undiscounted Gittins index, 
we next discuss it briefly.
Denote by $\nu_i^{\beta, *}$ the discounted Gittins index for
 discount factor $0 < \beta < 1$.
 \citet{kellymab} showed that $\nu_i^{\beta, *}$ is nondecreasing 
in $\beta$.
From  (\ref{eq:girep}), it follows that the
$\nu_i^{\beta,  *}$'s are uniformly bounded above by any upper bound
on active rewards. 
Hence, there exists a unique finite limiting index obtained as $\beta$
 approaches one:
$\nu^*_i \triangleq \lim_{\beta \nearrow 1} \nu_i^{\beta, *}$.

It is thus natural to consider the corresponding Taylor-series expansion 
\[
\nu_i^{\beta, *} = \nu^*_i + \gamma_i^* (1-\beta) + o(1-\beta) \quad 
\text{as } \beta \nearrow 1,
\]
whose validity was established by
\citet{karoth96}, who further showed
 that the limiting 
first- and second-order indices $\nu_i^*$ and $\gamma_i^*$ yield an
optimal policy for the multiarmed bandit problem under the 
\emph{average criterion}: highest priority is  awarded 
to a bandit with largest first-order index; ties are broken using the
second-order index.
\citet{kellymab} gave 
an earlier result in such vein.

It is thus of interest to compute 
$\nu^*_i$. 
Yet, if one tries to take limits as $\beta$ tends to one in 
the VWB, CP, or SE algorithms, they all break down, as
they involve divisions by $1-\beta$. 
In contrast, our 
FP(0) algorithm readily
computes the undiscounted Gittins index: one need simply set 
 $\beta = 1$ in the stated computations.

\section{Computational Study}
\label{s:ce}
It is well known that, in contemporary computers, the
 arithmetic-operation count of an algorithm need not be the prime driver of
its 
runtime performance.
Instead,  
the exponentially widening gap  since the 1990s between processor and
memory performance makes computation bottlenecks increasingly due to
memory-access times.
For a relevant discussion of such issues see, e.g., \citet{dongeijk}.
It is thus necessary to test whether the improved complexity
 of our FP algorithm translates into improved
runtimes.

We thus conducted a computational study measuring running
times on random instances of different sizes, using the
author's MATLAB implementations 
 of the algorithms described herein.
The experiments were run on an HP xw9300 2.8 GHz AMD Opteron
workstation with 4 GB of RAM, using MATLAB R2006a 64-bit on Windows xp x64.

 The results are reported in 
Table \ref{tab:cr}. 
For each of the stated state space sizes $n$, a random problem 
instance was generated, 
and the elapsed times $t^{\textup{FP(0)}}$, $t^{\textup{FP(1)}}$, $t^{\textup{CP}}$,
and $t^{\textup{SE}}$ (in seconds)
expended to compute the 
$n$ Gittins index values were
recorded for the FP(0), FP(1), CP and 
SE algorithms, respectively. 
The table further shows the speedup factors
$t^{\textup{FP(1)}}/t^{\textup{FP(0)}}$, 
$t^{\textup{CP}}/t^{\textup{FP(0)}}$, 
$t^{\textup{SE}}/t^{\textup{FP(0)}}$, and $t^{\textup{FP(1)}}/t^{\textup{CP}}$.

The results
 show that the theoretical speedup factors of $2$ and $3/2$ of the 
FP(0) algorithm against the FP(1) and SE algorithms, respectively, slightly underestimate the 
measured factors.
They further show that the theoretical speedup factor of $3$ of
the FP(0) algorithm relative to the CP
algorithm strongly overestimates the measured factor, which lies 
sligthly above that against the SE
algorithm.
Further,  despite the FP(1) algorithm's
smaller arithmetic-operation count relative to the CP algorithm
 ($(4/3) n^3$ vs.\ $2 n^3$), it is strongly outperformed
 by the latter. 
Overall, the FP(0) algorithm is the clear winner in such an experiment,
 followed by the SE algorithm, which shows only slight improvements
 over the CP algorithm.

The observed discrepancy
 between theoretical  and measured speedup
factors is due to memory-access issues. 
The computation bottleneck of the
CP algorithm, as revealed by profiling, is due to the update of
matrix $\mathbf{A}^{(k)}$.
In contrast, 
the FP(0), FP(1), and SE algorithms have two computation bottlenecks
each, corresponding to the major matrix updates at each
step. 
In our implementation, 
the latter are performed on indexed submatrices of preallocated
 matrices, involving expensive
 noncontiguous, random-stride memory-access patterns.
In contrast, matrix $\mathbf{A}^{(k)}$ is efficiently
accessed as a contiguous memory block.

\begin{table}
\caption{Runtime (secs.) Comparison of Index Algorithms.}
\begin{center}
\begin{tabular}{ccrcrcrcrcccc} \hline\noalign{\smallskip}
$n$ & & $t^{\textup{FP(0)}}$ & & $t^{\textup{FP(1)}}$ && $t^{\textup{CP}}$
&  & $t^{\textup{SE}}$ & $\displaystyle \frac{t^{\textup{FP(1)}}}{t^{\textup{FP(0)}}}$ & 
$\displaystyle \frac{t^{\textup{CP}}}{t^{\textup{FP(0)}}}$ &
$\displaystyle \frac{t^{\textup{SE}}}{t^{\textup{FP(0)}}}$ & 
$\displaystyle \frac{t^{\textup{FP(1)}}}{t^{\textup{CP}}}$\\ 
\noalign{\smallskip} \hline \noalign{\smallskip}
$1000$ & & $15.3$  & & $33.0$ & & $24.6$ & & $24.1$ & $2.16$ &  $1.61$
& $1.58$ & $1.34$ \\
$1500$ & & $49.3$  & & $113.3$ & & $82.3$ & & $80.6$ & $2.30$ &
$1.67$ & $1.63$ & $1.38$ \\
$2000$ & & $118.1$ & & $270.1$ & &  $195.8$ & & $194.6$ & $2.29$ &
$1.66$ & $1.65$ & $1.38$\\
$2500$ & & $230.1$ & & $520.6$ & &  $379.9$ & & $377.3$ & $2.26$ &
$1.65$ & $1.64$ & $1.37$ \\
$3000$ & & $395.7$ & & $914.0$ & &  $659.6$ & & $648.6$ & $2.31$ &
$1.68$ & $1.66$ & $1.39$ \\
$3500$ & & $627.4$ & & $1429.0$ & &  $1051.6$ & & $1043.7$& $2.28$ &
$1.68$ & $1.66$ & $1.36$ \\
$4000$ & & $937.6$ & & $2122.8$ & &  $1568.6$ & & $1559.7$& $2.26$ &
$1.67$ & $1.66$ & $1.35$ \\ 
$4500$ & & $1347.0$& & $3064.6$ & &  $2250.5$ & & $2217.0$& $2.28$&
$1.67$ & $1.65$ & $1.36$ \\ 
$5000$ & & $1832.7$& & $4225.6$ & &  $3092.6$ & & $3051.5$& $2.31$ &
$1.68$ & $1.67$ & $1.37$ \\ 
$5500$ & & $2458.2$& & $5629.6$ & &  $4136.9$ & & $4058.7$& $2.29$&
$1.68$ & $1.65$ & $1.36$ \\
$6000$ & & $3195.3$& & $7283.8$ & &  $5402.2$ & & $5373.6$& $2.28$ &
$1.69$ & $1.68$ & $1.35$ \\
\hline
\end{tabular}
\end{center}
\label{tab:cr}
\end{table}

\section{Conclusions}
\label{sec-Conclusions}
We have used the parametric version of Dantzig's
simplex method as the basis for designing a new algorithm to
compute the Gittins index of an $n$-state bandit, having an improved
operation count of $(2/3) n^3$, which
matches that of
solving a linear-equation system by Gaussian
elimination.
The algorithm has further been shown to outperform 
alternative methods. 
However,  direct implementation using
indexed submatrices results in expensive noncontiguous,
random-stride memory-access patterns, which represent its computation bottleneck.
It would be worth investigating whether the latter can be 
reduced 
through advanced approaches in numerical linear algebra, such
as block-partitioned implementations that
exploit advanced-architecture computers, as discussed in
 \citet{dongeijk}.
Corresponding simplex-based algorithms are developed in 
\citet{nmsc06} for restless bandits, which can change state when
passive.
We have further introduced a novel Gittins-index solution 
to the classical problem of optimal stopping of a
Markov chain, which renders the new algorithm applicable
to the latter.
Such an approach yields new tools for the study of optimal-stopping
problems that, we believe, warrant further
investigation.

\section*{Acknowledgments}
The author thanks the anonymous Associate Editor and two reviewers for
valuable suggestions that led to improvements in the paper.
This work was
supported by the Spanish Ministry of Education \& Science
under a Ram\'on y Cajal Investigator Award and 
grant MTM2004-02334, by the European Commission's Network of Excellence 
Euro-NGI, and by the Autonomous Community of Madrid-UC3M through grant
UC3M-MTM-05-075.

\bibliographystyle{ijocv081}

\end{document}